\documentclass[12pt,reqno]{amsart}
\usepackage{enumerate}
\usepackage{hyperref}
\usepackage{fullpage} 

\usepackage[english]{babel}

\newenvironment{Proof}[1][Proof]{\par\noindent\textbf{#1.}~}
  {\hfill$\square$\smallskip\par}
\newcommand{\dx}{\mathrm{d}}

\newcommand{\E}{\mathcal{E}}

\newcommand{\singseries}{\mathfrak{S}}
\newcommand{\Odip}[2]{\mathcal{O}_{#1}\!\left(#2\right)\mathchoice{\!}{}{}{}}

\newcommand{\Odi}[1]{\Odip{}{#1}}

\newtheorem{Theorem}{Theorem}
\newtheorem{Lemma}{Lemma}

\newtheoremstyle{Nonumtheorems}
{10pt}
{6pt}
{\itshape}
{}
{\bfseries}
{.}
{.5em}
{\thmname{#1}\thmnote{ (#3)}}

\theoremstyle{Nonumtheorems}
\newtheorem{Nonumlemma}{Lemma}

\allowdisplaybreaks

\title[Goldbach representations of an integer]{The number of Goldbach representations of an integer}
\author[Languasco \lowercase{and}  Zaccagnini]{Alessandro Languasco \lowercase{and} Alessandro Zaccagnini}

\begin{document}
\selectlanguage{english}
\maketitle
\section{Introduction}

Let $\Lambda$ be the von Mangoldt function and
\[
R(n) = \sum_{h_1+h_2=n} \Lambda(h_1)\Lambda(h_2)
\]
be the counting function for the Goldbach numbers.
This paper is devoted to study the behaviour of the average
order of magnitude of $R(n)$ for $n\in [1,N]$, where $N$ is a large
integer.
We have the following
\begin{Theorem}
\label{Main-Th}
Let $N \geq 2$ and assume  the Riemann Hypothesis (RH) holds.
Then
\[
\sum_{n=1}^{N}
R(n)
=
\frac{N^{2}}{2}
-2
\sum_{\rho} \frac{N^{\rho + 1}}{\rho (\rho + 1)}
+
\Odi{N \log^{3}N
},
\]
where $\rho=1/2+i\gamma$ runs over
the non-trivial zeros of the Riemann zeta function $\zeta(s)$.
\end{Theorem}

The first result of this kind was proved in 1991 by Fujii who
subsequently improved it (see \cite{Fujii1991}-\cite{Fujii1991a}-\cite{Fujii1991b})
until reaching the error term
$\Odi{(N \log N)^{4/3}}$.
Then Granville \cite{Granville2007}-\cite{Granville2008b}
gave an alternative proof of the same result and, finally,
Bhowmik and Schlage-Puchta \cite{BhowmikS2010a}
were able to reach the error term
$\Odi{N \log^{5} N}$; in \cite{BhowmikS2010a}
they also proved that the error term is  $\Omega(N\log \log N)$.

Our result improves the upper bound in Bhowmik and Schlage-Puchta \cite{BhowmikS2010a} by a factor $\log^{2}N$.
In fact, this seems to be the limit of the method in the current state
of the circle-method technology: see the remark after the proof.

If one admits the presence of some suitable weight in our average,
this loss can be avoided.
For example, using the Fej\'er weight we could work with
$L(N;\alpha)=\sum_{n=-N}^N (N-\vert n \vert) e(n\alpha) =\vert T(N;\alpha) \vert^2$
instead of $T(N;\alpha)$ in \eqref{circle}. The key property is that, for
$1/N<\vert \alpha \vert \leq 1/2$,  the function
$L(N;\alpha)$ decays as $\alpha^{-2}$ instead of $\vert \alpha \vert^{-1}$
and so the dissection argument in \eqref{I3-estim} is now more efficient and does not cause any loss
of logs.
Such a phenomenon is well-known from the literature
about the existence of Goldbach numbers in short intervals,
see, \emph{e.g.}, Languasco and Perelli \cite{LanguascoP1994}.

In fact  we will obtain Theorem \ref{Main-Th} as a consequence
of a weighted result. Letting $\psi(x)=\sum_{m\leq x} \Lambda(m)$, we have
\begin{Theorem}
\label{average-Th}
Let $2\leq y \leq N$ and assume  the Riemann Hypothesis (RH) holds.
Then
\begin{equation}
\label{average}
\max_{y \in [2, N]}
\left\vert
\sum_{n=1}^{y}
\Bigl[
R(n) - (2\psi(n)-n)
\Bigr]
e^{-n/N}
\right\vert
\ll
N\log^3 N.
\end{equation}
\end{Theorem}

The key reason why we are able to derive Theorem \ref{Main-Th}
from Theorem \ref{average-Th} via partial summation
is that the exponential weight in \eqref{average}
just varies in the range $[e^{-1/N}, e^{-1}]$
and so it does not change the order of magnitude
of the functions involved.

We will  use the original Hardy and Littlewood \cite{HardyL1923} circle
method setting, \emph{i.e.},
the weighted exponential sum
\begin{equation}
\label{tildeS-def}
\widetilde{S}(\alpha)
=
\sum_{n=1}^{\infty}
\Lambda(n) e^{-n/N}
e(n\alpha),
\end{equation}
where $e(x)=\exp(2\pi i x)$,
since it lets us avoid the use of Gallagher's Lemma
(Lemma 1 of \cite{Gallagher1970}) and hence,
in  this conditional case,  it gives slightly sharper results, see Lemma \ref{LP-Lemma} below.
Such a function was also used by Linnik \cite{Linnik1946,Linnik1952}.
The new ingredient in this paper is Lemma \ref{I2-lemma} below
in which we unconditionally detect the existence of the term
$-2\sum_{\rho} N^{\rho + 1}/(\rho (\rho + 1))$ by
solving an arithmetic problem connected with the original
one (see eq.~\eqref{cancellation} below).
In the previously mentioned papers this is obtained applying
the explicit formula for $\psi(n)$ twice.

The ideas that lead to Theorem \ref{Main-Th} and \ref{average-Th}
work also for the sum of $k\geq 3$ primes, \emph{i.e.},
for the function
\[
R_k(n) = \sum_{h_1+\dotsc+h_k=n}
\Lambda(h_1)\dotsm\Lambda(h_k).
\]
We can prove the following
\begin{Theorem}
\label{Main-Th-k-primes}
Let  $k \geq 3$ be an integer, $N \geq 2$ and assume  the Riemann Hypothesis (RH) holds.
Then
\[
\sum_{n=1}^{N}
R_k(n)
=
\frac{N^{k}}{k!}
-k
\sum_{\rho} \frac{N^{\rho + k-1}}{\rho (\rho + 1)\dotsm (\rho + k - 1)}
+
\Odip{k}{N^{k-1} \log^{k}N
},
\]
where $\rho=1/2+i\gamma$ runs over
the non-trivial zeros of $\zeta(s)$.
\end{Theorem}
The proof of Theorem \ref{Main-Th-k-primes}
is completely similar to the one of Theorems \ref{Main-Th} and \ref{average-Th}.
We just remark that the main differences are in the use of the explicit formula for
\[
\psi_{j}(t)
:=
\frac{1}{j!}
\sum_{n\leq t}
(t-n)^{j}
\Lambda(n)
\]
where $j$ is a non-negative integer, and of the following version of
Lemma 1 of \cite{Languasco2000a}:
\begin{Nonumlemma}
Assume  the Riemann Hypothesis (RH) holds.
Let $N \geq 2$, $z= 1/N-2\pi i\alpha$ and $\alpha\in [-1/2,1/2]$.
Then
\[
\Bigl\vert
\widetilde{S}(\alpha) -
\frac{1}{z}
\Bigr\vert
\ll
N^{1/2}
\Bigl(
1+ (N \vert \alpha \vert)^{1/2}
\Bigr)
\log N .
\]
\end{Nonumlemma}

Another connected problem we can address with this technique is
a short-interval version of Theorem \ref{Main-Th}.
We can prove the following
\begin{Theorem}
\label{short-Th}
Let $2\leq H\leq N$ and assume  the Riemann Hypothesis (RH) holds.
Then
\[
\sum_{n=N}^{N+H}
R(n)
=
HN + \frac{H^2}{2}
-
2
\sum_{\rho} \frac{(N+H)^{\rho + 1}-N^{\rho + 1}}{\rho (\rho + 1)}
+
\Odi{N\log^{2}N \log H},
\]
where $\rho=1/2+i\gamma$ runs over
the non-trivial zeros of $\zeta(s)$.
\end{Theorem}

Also in this case we do not give a proof of Theorem \ref{short-Th};
we just remark that the main difference is in the use of the
exponential sum
$\sum_{n=N}^{N+H} e(n\alpha)$
instead of $\sum_{n=1}^{N} e(n\alpha)$.

\medskip
\textbf {Acknowledgments.}
We would like to thank Alberto Perelli for a discussion.

\section{Setting of the circle method}

For brevity, throughout the paper we write
\begin{equation}
\label{def-z}
  z
  =
  \frac1N - 2 \pi i \alpha
  \qquad
  \text{where $N$ is a large integer and }
  \alpha \in [-1/2, 1/2].
\end{equation}

The first lemma is a $L^{2}$-estimate for the
difference  $\widetilde{S}(\alpha) - 1/z$.
\begin{Lemma}[Languasco and Perelli \cite{LanguascoP1994}]
\label{LP-Lemma}
Assume RH. Let $N$ be a
sufficiently large integer and $z$ be as in \eqref{def-z}.
For  $0 \leq \xi \leq 1/2$, we have
\[
\int_{-\xi}^{\xi}
\Bigl\vert
\widetilde{S}(\alpha) - \frac{1}{z}
\Bigr\vert^{2}
\dx \alpha
\ll
N\xi \log^{2} N.
\]
\end{Lemma}

This follows immediately from the proof of
Theorem 1 of \cite{LanguascoP1994} since
the quantity we would like to estimate here is
$\widetilde{R}_{1}+\widetilde{R}_{3}+\widetilde{R}_{5}$ there.

Lemma \ref{LP-Lemma} is the main reason why we
use $\widetilde{S}(\alpha)$
instead of its truncated form
$S(\alpha) =
\sum_{n=1}^N
\Lambda(n)
e(n\alpha)
$
as in Bhowmik and Schlage-Puchta \cite{BhowmikS2010a}.
In fact Lemma \ref{LP-Lemma}
lets us avoid the use of Gallagher's
Lemma \cite{Gallagher1970} which leads to a loss
of a factor $\log^2 N$ in the final estimate
(compare Lemma \ref{LP-Lemma} with Lemma 4 of  \cite{BhowmikS2010a}).
For a similar phenomenon in a slightly
different situation see also Languasco \cite{Languasco2000a}.

The next four lemmas do not depend on RH.
By the residue theorem one can obtain
\begin{Lemma}[Eq.~(29) of \cite{LanguascoP1994}]
\label{residue}
Let $N \geq 2$, $1\leq n \leq N$ and $z$ be as in \eqref{def-z}.
We have
\[
\int_{-\frac{1}{2}}^{\frac{1}{2}}
\frac{e(-n\alpha)}{z^2}
\ \dx \alpha
=
ne^{-n/N} + \Odi{1}
\]
uniformly for every $n \leq N$.
\end{Lemma}

\begin{Lemma}
\label{incond-mean-square}
Let $N$ be a sufficiently large integer and $z$ be as in \eqref{def-z}.
We have
\[
\int_{-\frac{1}{2}}^{\frac{1}{2}}
\Bigl\vert
\widetilde{S}(\alpha) - \frac{1}{z}
\Bigr\vert^{2}
\ \dx \alpha
=
\frac{N}{2} \log N
+
\Odi{N (\log N)^{1/2}}.
\]
\end{Lemma}
\begin{Proof}
By the Parseval theorem and the Prime Number Theorem we have
\[
\int_{-\frac{1}{2}}^{\frac{1}{2}}
\vert
\widetilde{S}(\alpha)
\vert^{2}
\ \dx \alpha
=
 \sum_{m = 1}^{\infty} \Lambda^{2}(m) e^{-2m / N}
=
\frac{N}{2} \log N
+
\Odi{N}.
\]
Recalling that the equation at the beginning of page 318
of \cite{LanguascoP1994} implies
\[
\int_{-\frac{1}{2}}^{\frac{1}{2}}
\frac{\dx \alpha}{\vert z\vert^{2}}
=
\frac{N}{\pi} \arctan(\pi N),
\]
the Lemma immediately follows using the relation
$\vert a - b \vert^{2} = \vert a\vert^{2} + \vert b \vert^{2} - 2\Re(a\overline{b})$
and the Cauchy-Schwarz inequality.
\end{Proof}

Let
\begin{equation}
\label{V-def}
  V(\alpha)
  =
  \sum_{m = 1}^{\infty} e^{-m / N} e(m \alpha)
  =
  \sum_{m = 1}^{\infty} e^{-m z}
  =
  \frac1{e^z - 1}.
\end{equation}

\begin{Lemma}
\label{V-behaviour}
If $z$ satisfies \eqref{def-z} then $V(\alpha) = z^{-1} + \Odi{1}$.
\end{Lemma}

\begin{Proof}
We recall the that the function $w / (e^w - 1)$ has a power-series
expansion with radius of convergence $2 \pi$ (see for example Apostol
\cite{Apostol1976}, page 264).
In particular, uniformly for $|w| \le 4 < 2 \pi$ we have
$w / (e^w - 1) = 1 + \Odi{|w|}$.
Since $z$ satisfies \eqref{def-z} we have $|z| \le 4$ and the result
follows.
\end{Proof}

Let now
\begin{equation}
\label{T-def}
T(y;\alpha) =
\sum_{n=1}^{y}
e(n\alpha)
\ll \min\left(y; \frac{1}{\Vert \alpha\Vert}\right).
\end{equation}

\begin{Lemma}
\label{I2-lemma}
Let $N$ be a large integer, $2\leq y \leq N$ and $z$ be as in
\eqref{def-z}.
We have
\begin{align}
\label{I2-eval}
  \int_{-1/2}^{1/2} T(y; -\alpha)
\frac{(\widetilde{S}(\alpha)-1/z)}z \, \dx \alpha
=
\sum_{n=1}^{y} e^{-n / N} (\psi(n) - n)
  +
  \Odi{(y N \log N)^{1/2}}.
\end{align}
\end{Lemma}
We remark that Lemma \ref{I2-lemma} is unconditional
and hence it implies, using also Lemma \ref{sum-integral}, that the
ability of detecting the term depending on the zeros
of the Riemann $\zeta$-function
in Theorem \ref{Main-Th} does not depend on RH.

\begin{Proof}
Writing
$\widetilde{R}(\alpha)=\widetilde{S}(\alpha) - 1/z$,
by Lemma \ref{V-behaviour} we have
\begin{align}
\notag
  \int_{-1/2}^{1/2}  T(y; -\alpha)  \frac{\widetilde{R}(\alpha)}z \, \dx \alpha
  & =
  \int_{-1/2}^{1/2} T(y; -\alpha) \widetilde{R}(\alpha) V(\alpha) \, \dx \alpha
+
  \Odi{\int_{-1/2}^{1/2} |T(y; -\alpha)| \, |\widetilde{R}(\alpha)| \, \dx \alpha}
  \\
    & =
  \int_{-1/2}^{1/2} T(y; -\alpha) \widetilde{R}(\alpha) V(\alpha) \, \dx \alpha
\label{first-step}  +
  \Odi{  (y N \log N)^{1/2}},
\end{align}
since, by the Parseval theorem and Lemma \ref{incond-mean-square},
the error term above is
\[
  \ll
  \Bigl( \int_{-1/2}^{1/2} |T(y; -\alpha)|^2 \, \dx \alpha \Bigr)^{1/2}
  \Bigl( \int_{-1/2}^{1/2} |\widetilde{R}(\alpha)|^2 \, \dx \alpha \Bigr)^{1/2}
  \ll
  (y N \log N)^{1/2}.
\]
Again by Lemma \ref{V-behaviour}, we have
\[
  \widetilde{R}(\alpha)
  =
  \widetilde{S}(\alpha)
  -
  \frac1z
  =
   \widetilde{S}(\alpha)
  -
  V(\alpha)
  +
  \Odi{1}
\]
and hence \eqref{first-step} implies
\begin{align}
\notag
  \int_{-1/2}^{1/2}  T(y; -\alpha) \frac{\widetilde{R}(\alpha)}z \, \dx \alpha
& =
  \int_{-1/2}^{1/2}
    T(y; -\alpha)
    \bigl(\widetilde{S}(\alpha) - V(\alpha) \bigr) V(\alpha) \, \dx \alpha \\
\label{second-step}
  &
\quad
+
  \Odi{\int_{-1/2}^{1/2} |T(y; -\alpha)| \,  |V(\alpha)| \, \dx \alpha}
  +
  \Odi{(y N \log N)^{1/2}}.
\end{align}
The Cauchy-Schwarz inequality and the Parseval theorem imply that
\begin{align}
\notag
  \int_{-1/2}^{1/2} |T(y; -\alpha)| \, |V(\alpha)| \, \dx \alpha
  &\le
  \Bigl( \int_{-1/2}^{1/2} |T(y; -\alpha)|^2 \, \dx \alpha \Bigr)^{1/2}
  \Bigl( \int_{-1/2}^{1/2} |V(\alpha)|^2 \, \dx \alpha \Bigr)^{1/2} \\
  &\ll
\label{err-second-step}
  \Bigl( y \sum_{m = 1}^{\infty} e^{-2 m / N} \Bigr)^{1/2}
  \ll
  (y N)^{1/2}.
\end{align}
By \eqref{second-step}-\eqref{err-second-step}, we have
\begin{equation}
\label{start-third-step}
  \int_{-1/2}^{1/2}  T(y; -\alpha) \frac{\widetilde{R}(\alpha)}z \, \dx \alpha
  =
  \int_{-1/2}^{1/2}
    T(y; -\alpha)
    \bigl(\widetilde{S}(\alpha) - V(\alpha) \bigr) V(\alpha) \, \dx \alpha
  +
  \Odi{(y N \log N)^{1/2}}.
\end{equation}
Now, by \eqref{tildeS-def} and \eqref{V-def}, we can write
\[
  \widetilde{S}(\alpha)
  -
  V(\alpha)
  =
  \sum_{m = 1}^{\infty} (\Lambda(m) - 1) e^{-m / N}  e(m \alpha)
\]
so that
\begin{align}
\notag
  &\int_{-1/2}^{1/2}
    T(y; -\alpha)
    \bigl(\widetilde{S}(\alpha) - V(\alpha) \bigr) V(\alpha) \, \dx \alpha \\
\notag
  &=
  \sum_{n = 1}^{y}
 \sum_{m_1 = 1}^{\infty} (\Lambda(m_1) - 1) e^{-m_1 / N}
   \sum_{m_2 = 1}^{\infty} e^{-m_2 / N}
   \int_{-1/2}^{1/2} e((m_1 + m_2 - n) \alpha) \, \dx \alpha \\
\notag
 &=
  \sum_{n = 1}^{y}
 \sum_{m_1 = 1}^{\infty} (\Lambda(m_1) - 1) e^{-m_1 / N}
   \sum_{m_2 = 1}^{\infty} e^{-m_2 / N}
   \begin{cases}
     1 & \text{if $m_1 + m_2 = n$} \\
     0 & \text{otherwise}
   \end{cases} \\
\label{cancellation}
  &=
  \sum_{n = 1}^{y}
    e^{-n / N}
    \sum_{m_1 = 1}^{n - 1} (\Lambda(m_1) - 1)
  =
  \sum_{n = 1}^{y}
    e^{-n / N} (\psi(n - 1) - (n - 1)),
\end{align}
since the condition $m_1 + m_2 = n$ implies that both variables are $< n$.
Now $\psi(n) = \psi(n - 1) + \Lambda(n)$, so that
\[
  \sum_{n = 1}^{y} e^{-n / N} (\psi(n - 1) - (n - 1))
  =
  \sum_{n = 1}^{y} e^{-n / N} (\psi(n) - n)
  +
  \Odi{y}.
\]
By \eqref{start-third-step}-\eqref{cancellation} and the previous
equation, we have
\begin{align*}
  \int_{-1/2}^{1/2}  T(y; -\alpha) \frac{\widetilde{R}(\alpha)}z \, \dx \alpha
  &=
  \sum_{n = 1}^{y} e^{-n / N} (\psi(n) - n)
  +
  \Odi{y + (y N \log N)^{1/2}} \\
  &=
  \sum_{n = 1}^{y} e^{-n / N} (\psi(n) - n)
  +
  \Odi{(y N \log N)^{1/2}}
\end{align*}
since $y \le N$, and hence \eqref{I2-eval} is proved.
\end{Proof}

\begin{Lemma}
\label{sum-integral}
Let $M > 1$ be an integer. We have that
\[
	\sum_{n = 1}^M (\psi(n)-n)
	=
	-
	\sum_{\rho} \frac{M^{\rho + 1}}{\rho (\rho + 1)}
	+
 	\Odi{M}.
\]
\end{Lemma}
\begin{Proof}
We recall the definition of $\psi_{0}(t)$ as $\psi(t)-\Lambda(t)/2$ if $t$ is an integer
and as $\psi(t)$ otherwise. Hence
\[
	\sum_{n = 1}^M \psi(n)
	=
	\sum_{n = 1}^M \psi_{0}(n)
	+
	\frac{1}{2}\sum_{n = 1}^M \Lambda(n)
	=
	\sum_{n = 1}^M \psi_{0}(n)
	+
 	\Odi{M}
\]
by the Prime Number Theorem.
Using the fact that $\psi_{0}(n) = \psi_{0}(t)$ for every
$t\in(n,n+1)$, we also get
\[
	\sum_{n = 1}^M \psi_{0}(n)
	=
	\sum_{n = 1}^M \int_n^{n+1} \psi_{0}(t)\ \dx t
	=
	\int_{0}^{M}  \psi_{0}(t)\ \dx t
	+
 	\Odi{M}.
\]
Remarking that
\[
	\sum_{n = 1}^M n
	=
	\int_{0}^{M} t\ \dx t
	+
 	\Odi{M},
\]
we can write
\begin{equation}
\label{discreto-continuo}
	\sum_{n = 1}^M (\psi(n)-n)
	=
	\int_{0}^{M}  (\psi_{0}(t)-t)\ \dx t
	+
 	\Odi{M}
	=
	\int_{2}^{M}  (\psi_{0}(t)-t)\ \dx t
	+
 	\Odi{M},
\end{equation}
since the integral on $(0,2]$ gives a contribution $\Odi{1}$.
For $t\geq 2$ we will use the explicit formula
(see eq.~(9)-(10) of \S17 of Davenport \cite{Davenport2000})
\begin{equation}
\label{explicit-truncated}
\psi_{0}(t)
=
t - \sum_{\vert \gamma\vert  \le Z} \frac{t^\rho}{\rho}
-\frac{\zeta'}{\zeta}(0)
-\frac12 \log \left(1-\frac{1}{t^{2}}\right)
+
R_{\psi}(t,Z),
\end{equation}
where
\begin{equation}
\label{explicit-error}
R_{\psi}(t,Z)
\ll
\frac{t}{Z}\log^2(tZ)
+
(\log t) \min\left(1;\frac{t}{Z\Vert t \Vert}\right).
\end{equation}
The term $-\frac{\zeta'}{\zeta}(0)
-\frac12 \log \left(1-\frac{1}{t^{2}}\right)$ gives a contribution
$\Odi{M}$ to the integral over $[2,M]$ in \eqref{discreto-continuo}.
We need now a $L^1$ estimate of the error term defined in \eqref{explicit-error}.
Let
\begin{equation}
\label{err-def}
\E(M,Z)
:=
\int_2^{M}
\vert
R_{\psi}(t,Z)
\vert
\dx t.
\end{equation}
The first term in  \eqref{explicit-error} gives a total
contribution to $\E(M,Z)$
which is
\begin{equation}
\label{err1}
\ll
\frac{M^{2}}{Z} \log^{2}(MZ).
\end{equation}
The second term  in  \eqref{explicit-error} gives a total
contribution to $\E(M,Z)$ which is
\begin{align}
\notag
&
\ll
\log M
\sum_{n=2}^{M}
\Bigl[
\int_{n}^{n+1/2}
\min
\Bigl(
1; \frac{t}{Z(t-n)}
\Bigr)
\dx t+
\int_{n+1/2}^{n+1}
\min
\Bigl(
1; \frac{t}{Z(n+1-t)}
\Bigr)
\dx t
\Bigr]
\\
\notag
&
\ll
\log M
\sum_{n=2}^{M}
\Bigl(
\int_{n}^{n+1/Z}
\dx t
+
\int_{n+1/Z}^{n+1/2}
\frac{t \ \dx t}{Z(t-n)}
+
\int_{n+1-1/Z}^{n+1}
\dx t
+
\int_{n+1/2}^{n+1-1/Z}
\frac{t\ \dx t}{Z(n+1-t)}
\Bigr)
\\
&
\label{err2}
\ll
\log M
\sum_{n=2}^{M}
\Bigl(
\frac{3}{Z}
+
\frac{2n+1}{Z} \log(\frac{Z}{2})
\Bigr)
\ll
\frac{M^{2}}{Z} \log M \log Z.
\end{align}
Combining \eqref{err-def}-\eqref{err2}, for $Z= M\log^{2}M$ we have
that
\begin{equation}
\label{err-estim}
\E(M,Z)
\ll
M.
\end{equation}
Inserting now \eqref{explicit-truncated} and \eqref{err-estim} into
\eqref{discreto-continuo} we obtain
\begin{equation}
\label{finale-zeri-bassi}
	\sum_{n=1}^{M}  (\psi(n)-n)
	=
	-
	\sum_{\vert \gamma\vert  \le Z} \frac{M^{\rho + 1}}{\rho (\rho + 1)}
	+
 	\Odi{M}.
\end{equation}
The lemma follows from \eqref{finale-zeri-bassi},
by remarking that
\[
  \sum_{\vert \gamma\vert  > Z} \frac{M^{\rho + 1}}{\rho (\rho + 1)}
  \ll
  M^{2} \int_Z^{+\infty} \frac{\log t}{t^2} \, \dx t
  \ll
  \frac{M^{2} \log Z}Z
  \ll
    \frac{M}{\log M}
\]
since  $Z=M\log^{2} M$.
\end{Proof}

\section{Proof of Theorem \ref{Main-Th}}

We will get Theorem \ref{Main-Th} as a consequence of
Theorem \ref{average-Th}.
By partial summation we have
\begin{align}
\notag
\sum_{n=1}^{N}
\Bigl[
R(n) - (2\psi(n)-n)
\Bigr]
&=
\sum_{n=1}^{N}
e^{n/N}
\Bigl\{
\Bigl[
R(n) - (2\psi(n)-n)
\Bigr]
e^{-n/N}
\Bigr\}
\\
&
\notag
=
e
\sum_{n=1}^{N}
\Bigl[
R(n) - (2\psi(n)-n)
\Bigr]
e^{-n/N}
\\
&
\qquad
\label{psum}
-\frac{1}{N}
\int_0^N
\Bigl\{
\sum_{n=1}^{y}
\Bigl[
R(n) - (2\psi(n)-n)
\Bigr]
e^{-n/N}
\Bigr\}
e^{y/N}
\ \dx y
+\Odi{1}.
\end{align}
Inserting  \eqref{average} in \eqref{psum} we get
\begin{align*}
\sum_{n=1}^{N}
\Bigl[
R(n) - (2\psi(n)-n)
\Bigr]
\ll
N\log^3 N
\end{align*}
and hence
\begin{equation}
\label{BP-improved}
\sum_{n=1}^{N}
R(n)
=
\sum_{n=1}^{N}
n
+
2\sum_{n=1}^{N}
(\psi(n)-n)
+\Odi{
N\log^3 N
}.
\end{equation}
Theorem \ref{Main-Th} now follows inserting  Lemma \ref{sum-integral}
and the identity $\sum_{n=1}^{N}  n = N^{2}/2 + \Odi{N}$ in \eqref{BP-improved}.

\section{Proof of  Theorem \ref{average-Th}}
Let $2\leq y \leq N$.
We first recall the definition of the singular
series of the Goldbach problem:
$\singseries(k)=0$ for $k$ odd and
\[
  \singseries(k)
  =
  2 \prod_{p > 2} \Big( 1 - \frac{1}{(p-1)^2} \Big)
\prod_{\substack{p \mid k \\ p > 2}} \frac{p - 1}{p - 2}
\]
for $k$ even. Hence,
using the well known estimate
$R(n)\ll n {\mathfrak S}(n) \ll n \log \log n$,
we remark that
\begin{equation}
\label{stima-banale}
\sum_{n=1}^{y}
\Bigl[
R(n) - (2\psi(n)-n)
\Bigr]
e^{-n/N}
\ll
\sum_{n=1}^{y}
n \log \log n
\ll
y^2 \log \log y.
\end{equation}
So it is clear  that  \eqref{average} holds for every $y\in [2, N^{1/2}]$.

Assume now that $y\in [N^{1/2},N]$ and
let $\alpha\in[-1/2,1/2]$.
Writing
$\widetilde{R}(\alpha)=\widetilde{S}(\alpha) - 1/z$,
recalling \eqref{T-def}
we have
\begin{align}
\notag
\sum_{n=1}^{y}
e^{-n/N}R(n)
&=
\sum_{n=1}^{y}
\int_{-\frac{1}{2}}^{\frac{1}{2}}
\widetilde{S}(\alpha)^2
e(-n\alpha)
\ \dx \alpha
=
\int_{-\frac{1}{2}}^{\frac{1}{2}}
\widetilde{S}(\alpha)^2
T(y;-\alpha)
\ \dx \alpha
\\
\notag
&
=
\int_{-\frac{1}{2}}^{\frac{1}{2}}
\frac{T(y;-\alpha)}{z^2}
\ \dx \alpha
+
2
\int_{-\frac{1}{2}}^{\frac{1}{2}}
\frac{T(y;-\alpha)\widetilde{R}(\alpha)}{z}
\ \dx \alpha
+
\int_{-\frac{1}{2}}^{\frac{1}{2}}
T(y;-\alpha)\widetilde{R}(\alpha)^{2}
\ \dx \alpha
\\
\label{circle}
&
= I_{1}(y) + I_{2}(y) +I_{3}(y),
\end{align}
say.

\paragraph{\textbf{Evaluation of $I_{1}(y)$}}
By Lemma \ref{residue}
we obtain
\begin{align}
\notag
 I_{1}(y)
&
=
\int_{-\frac{1}{2}}^{\frac{1}{2}}
\frac{T(y;-\alpha)}{z^2}
\ \dx \alpha
=
\sum_{n=1}^{y}
\int_{-\frac{1}{2}}^{\frac{1}{2}}
\frac{e(-n\alpha)}{z^2}
\ \dx \alpha
=
\sum_{n=1}^{y}
\Bigl(
n e^{-n/N} +\Odi{1}
 \Bigr)
 \\
\label{I1-eval}
&=
\sum_{n=1}^{y}
n e^{-n/N}
+
\Odi{y}.
\end{align}

\paragraph{\textbf{Estimation of $I_{2}(y)$}}
By \eqref{I2-eval} of  Lemma \ref{I2-lemma}
we obtain
\begin{align}
\label{I2-final}
I_{2}(y)
=
2
\sum_{n=1}^{y}
e^{-n/N} (\psi(n)-n)
+
  \Odi{(y N \log N)^{1/2}}.
\end{align}

\paragraph{\textbf{Estimation of $I_{3}(y)$}}

Using \eqref{T-def} and Lemma \ref{LP-Lemma} we have that
\begin{align}
\notag
 I_{3}(y)
&
\ll
\int_{-\frac{1}{2}}^{\frac{1}{2}}
\vert T(y;-\alpha) \vert
\vert \widetilde{R}(\alpha)\vert^{2}
\ \dx \alpha
\ll
y
\int_{-\frac{1}{y}}^{\frac{1}{y}}
\vert \widetilde{R}(\alpha)\vert^{2}
\ \dx \alpha
+
\int_{\frac{1}{y}}^{\frac{1}{2}}
\frac{\vert \widetilde{R}(\alpha)\vert^{2}}{\alpha}
\ \dx \alpha
+
\int_{-\frac{1}{2}}^{-\frac{1}{y}}
\frac{\vert \widetilde{R}(\alpha)\vert^{2}}{\vert \alpha \vert}
\ \dx \alpha
\\
\notag
&
\ll
N \log^{2}N
+
\sum_{k=1}^{\Odi{\log y}}
\frac{y}{2^{k}}
\int_{\frac{2^{k}}{y}}^{\frac{2^{k+1}}{y}}
\vert \widetilde{R}(\alpha)\vert^{2}
\ \dx \alpha
\ll
N \log^{2}N
+
\sum_{k=1}^{\Odi{\log y}}
\frac{y}{2^{k}} N \frac{ 2^{k+1}}{y}  \log^{2} N
\\
\label{I3-estim}
&
\ll
N\log^{2}N  \log y
.
\end{align}

\paragraph{\textbf{End of the proof}}

Inserting \eqref{I1-eval} and \eqref{I2-final}-\eqref{I3-estim}
into \eqref{circle} we immediately have
\[
\sum_{n=1}^{y}
e^{-n/N}R(n)
=
\sum_{n=1}^{y}
n e^{-n/N}
+
2
\sum_{n=1}^{y}
e^{-n/N} (\psi(n)-n)
+
\Odi{
N\log^{2}N  \log y
}.
\]
Hence
\[
\sum_{n=1}^{y}
e^{-n/N}
\Bigl[
R(n) -(2\psi(n)-n)
\Bigr]
\ll
N\log^{2}N  \log y
\]
and the maximum of the right hand side is attained
at $y=N$.
Thus we can write
\begin{equation}
\label{stima-y-alti}
\max_{y\in [N^{1/2},N]}
\left \vert
\sum_{n=1}^{y}
\Bigl[
R(n) - (2\psi(n)-n)
\Bigr]
e^{-n/N}
\right \vert
\ll
N \log^{3}N.
\end{equation}
Combining \eqref{stima-banale} and \eqref{stima-y-alti}
we get that Theorem \ref{average-Th} is proved.

\paragraph{Remark}
Let
\[
  f(\alpha)
  =
  f_N(\alpha)
  =
  \begin{cases}
    \frac12 N^{1/2} \log N
      & \text{if $\Vert \alpha \Vert \le (\log N)^{-1}$,} \\
    0
      & \text{if $\Vert \alpha \Vert >   (\log N)^{-1}$.}
  \end{cases}
\]
Then $f$ satisfies both Lemma \ref{LP-Lemma} and Lemma
\ref{incond-mean-square} in the sense that
\[
  \int_{-\xi}^{\xi} |f(\alpha)|^2 \, \dx \alpha
  \ll
  N \xi (\log N)^2
\]
for all $\xi \in [0, 1/2]$, and
\[
  \int_{-1/2}^{1/2} |f(\alpha)|^2 \, \dx \alpha
  =
  \frac12 N \log N,
\]
but
\[
  \int_{1/y}^{1/2} \frac{|f(\alpha)|^2}{\alpha} \, \dx \alpha
  =
  \frac14
  N (\log N)^2
  \int_{1/y}^{1/\log N} \frac{\dx \alpha}{\alpha}
  =
  \frac14
  N (\log N)^2 \log(y / \log N)
  \asymp
  N (\log N)^3
\]
for $y = N^{1/2}$ and sufficiently large $N$.
This means that the crucial bound for $I_3(y)$ in \eqref{I3-estim} is
essentially optimal in the present state of knowledge, and that it can
not be improved without deeper information on
$\widetilde{S}(\alpha) - z^{-1}$, such as the stronger analogue of
Lemma~\ref{LP-Lemma} that follows from a suitable form of Montgomery's
Pair-Correlation Conjecture.

\vskip 1cm
\noindent
Alessandro Languasco, Dipartimento di Matematica Pura e Applicata, Universit\`a
di Padova, Via Trieste 63, 35121 Padova, Italy; languasco@math.unipd.it

\medskip
\noindent
Alessandro Zaccagnini, Dipartimento di Matematica, Universit\`a di Parma, Parco
Area delle Scienze 53/a, Campus Universitario, 43124 Parma, Italy;
alessandro.zaccagnini@unipr.it

\end{document}